\documentstyle[11pt]{article}

\title{Endomorphisms of Banach algebras of infinitely differentiable
functions on compact plane sets}
\author{Joel F. Feinstein and Herbert Kamowitz}

\begin{document}
\newcommand{\dis}{\displaystyle}
\newcommand{\pis}{\mbox{$\overline{\phi(b)}$}}
\def\unorm#1#2{{\|#1\|_{{#2}}}}
\def\norm#1{{\|#1\|}}
\def\rd{{\rm d}}

\maketitle


     This note is a sequel to \cite{nfact} where we investigated
the endomorphisms of a certain class of Banach algebras of 
infinitely differentiable functions on the unit interval.

    Start with a perfect, compact plane set $X$.
We say that a complex-valued function $f$ defined on $X$ 
is {\it complex-differentiable} 
at a point $a\in X$ if
the limit
$$ f^\prime(a) = \lim_{z\to a, \ z\in X}
{f(z)-f(a)\over z-a} $$
exists. We call $f^\prime(a)$ the {\it complex derivative} of $f$ at $a$.
Using this concept of derivative, we define the terms 
{\it complex--differentiable on} $X$, 
{\it continuously complex--differentiable on} $X$, 
and {\it infinitely complex--differentiable on} $X$
in the obvious way. We denote the $n$-th complex derivative
of $f$ at $a$ by $f^{(n)}(a)$, and we denote the set of infinitely 
differentiable functions on $X$ by $D^\infty(X)$.

    Let $M_n$ be a sequence of positive numbers satisfying $M_0=1$ and
$\dis \left(\begin{array}{c}m+n\\n \end{array} \right) \leq \frac{M_{m+n}}{M_mM_n}$, and let $\dis D(X,M)=\{f \in D^{\infty}(X): \|f\|=\sum_{n=0}^{\infty}
\frac{\|f^{(n)}\|_{\infty}}{M_n} < \infty\}.$ With pointwise addition
and multiplication, $D(X,M)$ is a commutative normed algebra which
is not necessarily complete.

     Clearly all polynomials
when restricted to $X$ belong to each $D(X,M)$. It was further proved
in \cite{d} that
the algebra $D(X,M)$ includes 
all the rational functions with poles off $X$ if and only if 
$$\lim_{n\to\infty} {\left({n!}\over{M_n}\right)}^{1\over n} = 0. \eqno(1)$$
We say that $M_n$ is a {\it nonanalytic sequence} if (1) holds \cite{d}.

   These algebras were considered by Dales and Davie in connection
with several then open questions. In \cite{d} examples of such Banach algebras
were constructed; one gave an example of a commutative semisimple
Banach algebra for which the peak points were of first category
in the Silov boundary, and a second  was an example of a commutative
semisimiple Banach algebra $B$ and a discontinuous function $F$
acting on $B.$

    From a different direction, we were led to these algebras
in connection with the study of compact endomorphisms of
commutative semisimple Banach algebras. On the basis of many
early examples, it was conjectured that all nonzero compact endomorphisms
of regular commutative semisimple Banach algebras with connected
maximal ideal spaces were essentially point evaluations.
Our first counterexample 
was the endomorphism $T:\dis Tf(x)=f(\frac{x}{2})$ on
$D([0,1],n!^2)$. \cite{hkcpt1}

     This led to the question of determining all
endomorphisms of $D([0,1],n!^2)$ and more generally 
determining the endomorphisms of other algebras 
$D(X,M)$. In \cite{nfact}, 
the question was settled for $X=[0,1]$ in many cases when the
weights $M_n$ were nonanalytic.



    In this paper we look at more general perfect, compact subsets of the plane, 
and 
investigate the extent to which the results for the interval extend to this
setting. In particular we shall give a variety of results in 
the case of the closed unit disk. We shall also partially resolve some of
the problems left open for the interval in \cite{nfact}.

      In general, the normed algebra $D(X,M)$ need not be complete. However, if
the compact set $X$ is a finite union of uniformly regular sets, 
\footnote{A compact plane set $X$ is {\it uniformly regular }
if, for
all $z, w \in X$, there is a rectifiable arc in $X$ joining $z$ to $w$,
and the metric given by the geodesic distance between points of $X$
is uniformly equivalent to the Euclidean metric. \cite{d}}
then $D(X,M)$ is complete
for every weight $M_n$. Such sets include $[0,1]$ and $\bar{\Delta}$
where $\Delta$ is the open unit disc.

    We recall further results from \cite{d}.  Suppose
$D_R(X,M)$ is the  closed subalgebra of $D(X,M)$
 generated by the rational
functions with poles off $X$.  If $M_n$ is nonanalytic, and
$X$ is uniformly regular, then $D_R(X,M)$ is {\it natural}
meaning that the maximal ideal space of $D_R(X,M)$ is $X$.
Further, for nicely shaped $X$, 
it was shown in \cite{flo} that $D_R(X,M)=D(X,M).$ 
More can be said in the cases of the unit disc, $\Delta$,
and unit interval.
Indeed, it was shown in \cite{dd} that the polynomials are
dense in $D(\bar{\Delta},M)$ and in \cite{of} that the polynomials
are dense in $D([0,1],M).$


In contrast, when the weight $M_n=n!$, the maximal ideal space
of the Banach algebra
$D([0,1],n!)$ 
equals $\{\lambda:{\rm dist}(\lambda,[0,1]) \leq 1\}.$
If $\dis \sum_{n=0}^{\infty} \frac{M_{n}}{M_{n+1}}=\infty$ then the
algebra is {\it quasi-analytic} in the sense that if $f^{(n)}(a)=0$ for
all $n$, then $f = 0$.

    Some examples when $X=[0,1]:$ 
If $\dis M_n=n!^{\alpha}$, $\alpha > 1,$
the algebra $D([0,1],M)$ is natural and regular, while  
if $\dis M_n=n![\log(n+1)]^n$ the algebra is natural and quasi-analytic. If $M_n=n!$
then the algebra is quasi-analytic and not natural.

    From general principles, if $T$ is a nonzero endomorphism
of a commutative semisimple Banach algebra $B$ with maximal ideal space
$\Phi_B$, 
then there exists a
continuous map $\phi:\Phi_B \rightarrow \Phi_B$ such that
$\widehat{Tf}(x)=\hat{f}(\phi(x))$ for all $f \in B$, $x \in \Phi_B.$
For a natural Banach algebra  $D(X,M)$, 
our question can be restated as
determining conditions on $\phi:X \rightarrow X$ for which
$f \circ \phi \in D(X,M)$ for all $f \in D(X,M).$


   The results in Parts A and B concern two types of theorems.
One will deal with arbitrary nonanalytic weights $M_n$ 
and $\phi:X \rightarrow X$ which satisfy 
an additional analyticity
property. Then the 
analyticity property on $\phi$ will be removed and weights
$M_n$ constructed such that the theorems hold. 
The paper will then conclude in Part C with a detailed study when
$X=\bar{\Delta}.$

   Finally, in the remainder of this note any unlabeled norm
$\|\cdot\|$ will denote the sup norm of the function on $X$.
The symbol $\|\cdot\|_{D(X,M)}$ will be used for the algebra
norm.

{\bf Part A}

   The following was proved in \cite{nfact} for $X=[0,1]$,
but an examination of the proof shows that the reasoning does not
depend on $[0,1]$. Since references will be made to this proof, 
for the convenience of the reader 
we reproduce the proof of part (a) from \cite{nfact} 
in the more general form.
We remark, too, that the proof of part (b) also goes through with minor
modifications.


{\bf Theorem 1:}  Let $X$ be a perfect, compact plane set
and $M_n$ be a nonanalytic weight sequence.

(a) Suppose $\phi \in D^{\infty}(X)$, $\phi:X \rightarrow X$, 
$\dis \limsup_{k \rightarrow \infty}
(\frac{\|\phi^{(k)}\|}{k!})^{1/k}$ is finite and $\|\phi'\|_{\infty} < 1.$ Then $\phi$ induces an endomorphism of $D(X,M)$.

(b) If, in addition, there exists $B>0$ such that $\dis \frac{M_m}{m!}\frac{n!}{M_n} \leq \frac{B}{m^{n-m}}$ for $n \geq m \geq 1$, and $\dis \{\frac{\|\phi^{(k)}\|}{k!}\}$ is bounded,
then $\|\phi'\|_{\infty} \leq 1$
is sufficient for $\phi$ to induce an endomorphism of $D(X,M)$.

{\bf Proof:}

(a) Suppose $\phi$ satisfies the hypotheses. Let $F \in D(X,M).$ We show that $F \circ \phi \in D(X,M).$  The 
following equality for higher derivatives of composite functions is known
as Fa\`{a} di Bruno's formula.

\[\dis \frac{d^n}{dx^n}(F \circ \phi)=\sum_{m=0}^{n} F^{(m)}(\phi)
{\bf \Sigma} \frac{n!}{a_1!a_2!\cdots a_n!} \frac{(\phi')^{a_1}(\phi'')^{a_2} \cdots (\phi^{(n)})^{a_n}}{1!^{a_1}2!^{a_2}\cdots n!^{a_n}}\]
where the inner sum {\bf ${\Sigma}$} is over non-negative integers $a_1, a_2, \cdots,
a_n$ satisfying $a_1+a_2+ \cdots +a_n=m$ and $a_1+2a_2+\cdots +na_n=n.$

{\bf Throughout the proof of Theorem 1 the inner sum $\Sigma$ will always
be over non-negative integers $a_1, a_2, \cdots,a_n$ satisfying $a_1+a_2+ \cdots +a_n=m$ and $a_1+2a_2+\cdots +na_n=n.$}

    Thus, Fa\`{a} di Bruno's formula implies that

\[\dis \|\frac{d^n}{dx^n}(F \circ \phi) \| \leq \sum_{m=0}^n \|F^{(m)}
(\phi)\| {\bf \Sigma}\frac{n!}{a_1!a_2!\cdots a_n!} \frac{(\|\phi'\|)^{a_1}(\|\phi''\|)^{a_2} \cdots (\|\phi^{(n)}\|)^{a_n}}{1!^{a_1}2!^{a_2}\cdots n!^{a_n}}\]
and so

\[\dis \sum_{n=0}^{\infty} \frac{1}{M_n}\|\frac{d^n}{dx^n}(F \circ \phi)\|
 \leq \sum_{n=0}^{\infty} \frac{1}{M_n} \sum_{m=0}^n \|F^{(m)}\|
 {\bf \Sigma}\frac{n!}{a_1!a_2!\cdots a_n!} \frac{(\|\phi'\|)^{a_1}(\|\phi''\|)^{a_2} \cdots (\|\phi^{(n)}\|)^{a_n}}{1!^{a_1}2!^{a_2}\cdots n!^{a_n}}.\]

Then, after interchanging the order of summation, we have

\[ \sum_{n=0}^{\infty} \frac{1}{M_n}\|\frac{d^n}{dx^n}(F \circ \phi)\|
 \leq \sum_{m=0}^{\infty} \|F^{(m)}\| \sum_{n=m}^{\infty} \frac{1}{M_n}
 {\bf \Sigma}\frac{n!}{a_1!a_2!\cdots a_n!} \frac{(\|\phi'\|)^{a_1}(\|\phi''\|)^{a_2} \cdots (\|\phi^{(n)}\|)^{a_n}}{1!^{a_1}2!^{a_2}\cdots n!^{a_n}},\]
whence


(*)

\[ \sum_{n=0}^{\infty} \frac{1}{M_n}\|\frac{d^n}{dx^n}(F \circ \phi)\|
 \leq \sum_{m=0}^{\infty} \frac{\|F^{(m)}\|}{M_m}m! \sum_{n=m}^{\infty} \frac{1}{n!} {\bf \Sigma}\frac{M_m}{m!} \frac{n!}{M_n}\frac{n!}{a_1!a_2!\cdots a_n!} \frac{(\|\phi'\|)^{a_1}(\|\phi''\|)^{a_2} \cdots (\|\phi^{(n)}\|)^{a_n}}{1!^{a_1}2!^{a_2}\cdots n!^{a_n}}.\]


    Choose $\epsilon > 0$ such that $\dis q:=\|\phi'\|(1+\sum_{k=2}^{\infty}\frac{\|\phi^{(k)}\|}{\|\phi'\|k!}\epsilon^{k-1}) < 1$. Such $\epsilon$ exists
since $\dis h(\lambda)=\sum_{k=2}^{\infty} \frac{\|\phi^{(k)}\|}{k!}\lambda^{k-1}$ is analytic near $0$ and $h(0)=0.$

     Since $\dis \lim_{k \rightarrow \infty} (\frac{k!}{M_k})^{1/k}=0$, 
there exists $B > 0$ such that $\dis \frac{k!}{M_k} < B\epsilon^k$ 
for all $k$, and since $\dis \frac{M_m}{m!}\frac{n!}{M_n} \leq \frac{(n-m)!}
{M_{n-m}}$ for $n \geq m$, we have that
$\dis (\frac{M_m}{m!})(\frac{n!}{M_n}) < B \epsilon^{n-m}$, $n \geq m.$

     Therefore the inequality (*) implies that

\[ \sum_{n=0}^{\infty} \frac{1}{M_n}\|\frac{d^n}{dx^n}(F \circ \phi)\|
 \leq \sum_{m=0}^{\infty} \frac{\|F^{(m)}\|}{M_m}m! \sum_{n=m}^{\infty} \frac{1}{n!} {\bf \Sigma} B\epsilon^{n-m}\frac{n!}{a_1!a_2!\cdots a_n!} \frac{\|\phi'\|^{a_1}\|\phi''\|^{a_2} \cdots \|\phi^{(n)}\|^{a_n}}{1!^{a_1}2!^{a_2}\cdots n!^{a_n}}\]

\[ = B \sum_{m=0}^{\infty} \frac{\|F^{(m)}\|}{M_m}m!\sum_{n=m}^{\infty}
\frac{\epsilon^n}{n!} {\bf \Sigma} \frac{n!}{a_1!a_2!\cdots a_n!} \frac{(\frac{\|\phi'\|}{\epsilon})^{a_1}(\frac{\|\phi''\|}{\epsilon})^{a_2} \cdots (\frac{\|\phi^{(n)}\|}{\epsilon})^{a_n}}{1!^{a_1}2!^{a_2}\cdots n!^{a_n}},\]
since $a_1+a_2+\cdots+a_n=m.$

     It follows from \cite{a}, page 823, B, formula 3, that

\[ m!\sum_{n=m}^{\infty}\frac{\epsilon^n}{n!} {\bf \Sigma}\frac{n!}{a_1!a_2!\cdots a_n!} \frac{(\frac{\|\phi'\|}{\epsilon})^{a_1}(\frac{\|\phi''\|}{\epsilon})^{a_2} \cdots (\frac{\|\phi^{(n)}\|}{\epsilon})^{a_n}}{1!^{a_1}2!^{a_2}\cdots n!^{a_n}}
=(\sum_{k=1}^{\infty} \frac{\|\phi^{(k)}\|}{\epsilon k!}\epsilon^k)^m.\]
     Moreover, from the definition of $q$,   
\[ (\sum_{k=1}^{\infty} \frac{\|\phi^{(k)}\|}{\epsilon k!}\epsilon^k)^m=
\|\phi'\|^m (1+ \sum_{k=2}^{\infty} \frac{\|\phi^{(k)}\|}{k!\|\phi '\|} \epsilon^{k-1})^m = q^m.\]

Therefore, since  $F \in D(X,M)$ and $0 < q < 1$, \[\sum_{n=0}^{\infty} \frac{1}{M_n}\|\frac{d^n}{dx^n}(F \circ \phi)\| \leq B \sum_{m=0}^{\infty} \frac{\|F^{(m)}\|}{M_m}q^m < \infty,\]
as required.



     The weights $M_n=n!^{\alpha}$, $\alpha \geq 2$ and $M_n=n!n^{n^2}$ are
nonanalytic weights for which the condition in (b),
$\dis \frac{M_m}{m!}\frac{n!}{M_n} \leq \frac{B}{m^{n-m}},$
holds. The condition fails for $M_n = n!^{\alpha}$, $1 < \alpha <2$,
\vspace{.005in}
and for $M_n=n!(\log(n+1))^n.$

     Throughout this paper, a self-map $\phi$  of $X$ will be said
to be {\it analytic on $X$} if $\phi$ extends to an analytic function
on a neighborhood of $X$.     

     The next theorem is an extension of Theorem 1(a) which is 
useful when $X$ has nonempty interior.

{\bf Theorem 2}
Let $X$ be a perfect, compact plane set. 
Let $\phi$ be an analytic
self-map of $X$.
 Let $M_n$ be a nonanalytic sequence.
Set $K=\{z\in X: |\phi'(z)|\geq 1 \}$, and
suppose that $\phi(K) \subset {\rm int}(X)$. 
Then $\phi$ induces an endomorphism of $D(X,M)$.

{\bf Proof:}

Choose $\epsilon > 0$ such that 
$\phi(\{z\in X:|\phi'(z)|\geq 1-\epsilon \}) \subset {\rm int}(X)$.
Set $K_1 = \{z\in X: |\phi'(z)| \leq 1-\epsilon\}$ 
and $K_2 = \{z\in X: |\phi'(z)| \geq 1-\epsilon\}$.
Then $X = K_1 \cup K_2$. For $f \in D(X,M)$, calculations similar
to those in the proof of Theorem 1(a) show
that the restriction of $f\circ\phi$ to
$K_1$ is in $D(K_1,M)$. Further $f \circ \phi$ is analytic on a 
neighborhood of $K_2$, whence the restriction of $f \circ \phi$ to
$K_2$ is in $D(K_2,M)$. It follows that $f \circ \phi \in D(X,M)$,
as required. 

   We note that the algebras in these theorems were not necessarily complete.

   This result may be improved if the stronger condition in part (b)
is placed on the sequence $M_n$. 
That is, suppose that $\phi$ is an analytic self-map of $X$ and $M_n$ 
satisfies the condition in part (b) of Theorem 1.
Then if 
$X = K_1 \cup K_2$
with $K_1$, $K_2$ compact such that $\phi(K_2) \subseteq {\rm int}(X)$
and $\unorm {\phi'} {K_2} \leq 1$, then the Fa\`a di Bruno 
calculations used in the proof of Theorem 1(b) in \cite{nfact},
combined with the argument above shows again that $\phi$
induces an endomorphism of $D(X,M)$.

     However, we shall see that  even in the case  where $X=\bar{\Delta},$
there are analytic self-maps $\phi$ on $\bar{\Delta}$
which have $|\phi'(z)| \leq 1$ for all those $z$ such that 
$|\phi(z)|=1$,
but for which no such decomposition into $K_1$ and $K_2$ is possible.

    In a different direction, if the infinitely differentiable 
self-maps $\phi$ are not 
as well behaved as in the previous hypothesis, 
the following result shows that 
maps $\phi$ still induce endomorphisms of $D(X,M)$
provided that the sequence $M_n$ grows rapidly enough. 

{\bf Theorem 3:}
Let $X$ be a perfect, compact plane set 
whose boundary is given by a finite union $\Gamma$ of 
piecewise smooth Jordan curves, and  such that
$\Gamma$ has winding number $1$ about each point of $X\backslash\Gamma$
and $0$ about each point of the complement of $X$.
Let $\phi \in D^\infty(X)$ with 
$\phi:X \rightarrow X$. Suppose that for all $z \in \phi^{-1}(\Gamma)$
we have $|\phi'(z)| \leq 1$. Then there exists a nonanalytic algebra 
sequence $M_n$ such that $\phi$ induces an endomorphism of $D(X,M)$.

{\bf Proof:}

   We choose the sequence $M_n$ inductively. We start with $M_0=M_1=1$.
For $n \geq 2$, and having chosen $M_0,M_1,\dots,M_{n-1}$, 
it follows from 
 Fa\`a di Bruno's formula that there are constants $C_{n,\phi,m} > 0$ 
such that for all infinitely differentiable functions $F$ on $X$ 
and all $z \in X$ we have 
$$|(F \circ \phi)^{(n)}(z)| \leq \|F^{(n)}\|_\infty|\phi'(z)|^n
    +\sum_{m=0}^{n-1}{C_{n,\phi,m}\|F^{(m)}\|_\infty}.$$

    Choose a relatively open set $U = U_{n,\phi} \supseteq \phi^{-1}(\Gamma)$ 
on which $|\phi'(z)|^n \leq 2$. Set $K=K_{n,\phi} = X \backslash U$.
Then $\phi(K)$ is a compact subset of the interior of $X$.
Assuming that $K_{n,\phi} \neq \emptyset$,
set $$A_{n,\phi} = \frac{L}{2\pi}
\sup\{|\frac{\rd^n}{\rd z^n}(\omega-\phi(z))^{-1}|:
   \omega \in \Gamma, z\in K_{n,\phi}\},$$
where $L$ is the total length of $\Gamma$.

    We see that, for $z \in K$ and $F$ as above, 
$$|(F \circ \phi)^{(n)}(z)|=
\left|{\frac{1}{2\pi i} \int_\Gamma{F(\omega) \frac{\rd^n}{\rd z^n}
(\omega-\phi(z))^{-1} \rd \omega}}\right|
\leq \|F\|_\infty A_{n,\phi}.$$

    Now choose $M_n$ large enough such that the following conditions 
are all satisfied:
(i) $M_n\geq (n!)^2$;
(ii) ${M_n\over{M_k M_{n-k}}} \geq {n \choose k}$ for $0 \leq k \leq n-1$;
(iii) $\sum_{m=0}^{n-1}{C_{n,\phi,m} M_m/M_n} \leq 2^{-n}$;
(iv) $A_{n,\phi}/M_n \leq 2^{-n}$.

     In the case where $K_{n,\phi} = \emptyset$, choose $M_n$ satisfying
(i) to (iii) above instead.

     Having chosen the sequence $M_n$ as above, we see that 
$M_n$ is clearly a nonanalytic 
sequence. For $F \in D(X,M)$ we have, for all $k\geq 0$,
$$\|F^{(k)}\|_{\infty} \leq  M_k \|F\|_{D(X,M)}.$$
Thus, by our choice of $M_n$, and considering separately the points $z$ in
$U_{n,\phi}$ and $K_{n,\phi}$, we have, for $n \geq 2$,


\[\frac{\|(F \circ \phi)^{(n)}\|_{\infty}}{M_n} \leq {\rm max}\{2^{-n}
\|F\|_{\infty},\frac{2\|F^{(n)}\|_{\infty}}{M_n}+2^{-n}\|F\|_{D(X,M)}\}\]
\[=\frac{2\|F^{(n)}\|_{\infty}}{M_n}+2^{-n}\|F\|_{D(X,M)}.\]
Thus $F\circ \phi$ is also in $D(X,M)$, as required. 

    An easy modification of this argument allows one sequence $M_n$ to
work for any given countable collection of such functions $\phi$ 
simultaneously. However these sequences $M_n$ will grow very rapidly.
The same argument can be used in an easier form to show that for
every self-map $\phi \in D^\infty([0,1])$ with $\|\phi'\|_{\infty} \leq 1$ there
is a nonanalytic sequence $M_n$, depending on $\phi$, such that $\phi$
induces an endomorphism of $D([0,1],M)$. Similarly we can show that
if the self-map $\phi \in D^\infty([0,1])$ with $\unorm {\phi'} \infty > 1$
there is a nonanalytic sequence $M$, depending on $\phi$, such that 
$\phi$ does not induce an endomorphism of $D([0,1],M)$. This will follow from
the more general result, Theorem 5.

    We now make a further definition. If $X$ is a subset of ${\bf C}$
and if $c \in X$, we say that $c$ has an {\it external circular tangent
} if there is an open disc $\Delta_1$ contained in the 
complement of $X$ with $\overline{\Delta_1} \cap X = \{c\}.$ It is easy to
see from the geometry that this condition is equivalent to each
of the following.

(i) For some $a \not\in X$, $|c-a|= {\rm min} \{|z-a|:z \in X\}$.

(ii) There exists a real number $\theta$ such that for large numbers $R$
we have $|1+e^{i\theta}R(c-z)| > 1$ for all $z \in X \backslash \{c\}.$

    Clearly every point of $[0,1]$ or  $\bar{\Delta} \setminus \Delta$
has an external circular tangent.

{\bf Part B}

    We now consider a converse to Theorem 1, and again state the
result for general perfect, compact subsets of ${\bf C}$. The
proof of this theorem is very similar to that of the corresponding
theorem in \cite{nfact}.

{\bf Theorem 4:} Let $M_n$ be a nonanalytic weight and suppose 
$\phi \in D^{\infty}(X)$,  $\phi:X 
\rightarrow X$ and $\dis \limsup_{n \rightarrow \infty}
(\frac{\|\phi^{(k)}\|}{k!})^{1/k}$ is finite.
Suppose for some $b \in X$, $\phi(b)$ has an external circular tangent
and $|\phi'(b)|> 1$. Also suppose that $D(X,M)$ is a 
Banach algebra.  Then $\phi$ does not induce an endomorphism
of $D(X,M).$

{\bf Proof:(outline)}

   Assume that $\phi$ induces an endomorphism, $|\phi'(b)| > 1$ 
and that $\phi(b)$ has an external
circular tangent. Let
   \[F_R(z)=\frac{1}{1+e^{i\theta}R(\phi(b)-z)}\]
where $\theta$ and $R$ are chosen so $\dis |1+e^{i\theta}R(\phi(b)-z)|>1$
for $z \in X \setminus \{\phi(b)\}$. Then $F_R \in D(X,M)$,
$\|F_R\|_{\infty}=1$
and $\dis \|F_R^{(m)}\|_{\infty} = m!R^m$. With a very slight 
modification, namely 
replacing $x_k$ by
$\rule{0pt}{6ex}\dis \frac{\overline{\phi'(b)}}{|\phi'(b)|}\frac{\phi^{(k)}(b)}{k!}$ 
rather than by $\dis \frac{\phi^{(k)}(b)}{k!}$, 
the proof proceeds
exactly 
as in  the proof of Theorem 3
of \cite{nfact}, eventually  arriving at \[\dis \|F_R\|_{D(X,M)}
=\sum_{m=0}^{\infty}\frac{m!R^m}{M_m}\] and for each $\epsilon$,
$0 < \epsilon < 1$, 
   \[\|F_R \circ \phi\|_{D(X,M)} \geq \sum_{m=0}^{\infty}\frac{((1-\epsilon)R|\phi'(b)|)^m}
{M_m}.\]

   Since  $M_n$ is a nonanalytic weight, $\dis \lim_{n \rightarrow \infty}
(\frac{n!}{M_n})^{1/n}=0$, and so  
$\dis g(z)=\sum_{n=0}^{\infty}\frac{n!}{M_n}z^n$
is a transcendental entire function. In general, if $\dis g(z)=
\sum_{n = 0}^{\infty} a_n z^n$ is a transcendental
entire function and if $\dis M_g(r)=\sup_{|z|=r}|g(z)|$, then for
$c > 1$, $\dis \lim_{r \rightarrow \infty} \frac{M_g(cr)}{M_g(r)} = \infty.$
(\cite{psz}, p 5, problem 24)

    Also since $D(X,M)$ is complete, every endomorphism is bounded.
Hence if $\phi$ induces an endomorphism, then there is a number $K > 0$
such that $\|f \circ \phi \|_{D(X,M)} \leq K \|f\|_{D(X,M)}$ for all $f \in D(X,M).$
Now for all $R > 0$, \[M_g(R)=g(R)=\|F_R\|_{D(X,M)}\]
and 
\[M_g((1-\epsilon)|\phi'(b)|R)=g((1-\epsilon)|\phi'(b)|R) \leq 
\|F_R \circ \phi\|_{D(X,M)}.\]


   Thus, if $\phi$ induces an endomorphism, then for some $K > 0$,
\[M_g((1- \epsilon)|\phi'(b)|R) \leq K M_g(R)\]
for large $R$, and so
\[\limsup_{R \rightarrow \infty} \frac{M_g((1-\epsilon)|\phi'(b)|R)}
{M_g(R)} \leq K < \infty.\]

    This implies that $\dis (1-\epsilon)|\phi'(b)| \leq 1.$ Letting
$\epsilon \rightarrow 0$ shows that $|\phi'(b)| \leq 1$ contrary
to our assumption.

    Remark: This theorem also has an interpretation when $D(X,M)$
is not complete. Suppose $\phi \in D^\infty(X)$,
$\phi:X \rightarrow X$ and $\dis \limsup_{n \rightarrow \infty}
(\frac{\|\phi^{(k)}\|}{k!})^{1/k}$ is finite.
Suppose for some $b \in X$, $\phi(b)$ has an external circular tangent
and $|\phi'(b)|> 1$. Then $\phi$ does not induce a {\it bounded}
endomorphism of $D(X,M).$

    For the next theorem, we call a compact plane set $X$ {\it good} if
$D(X,M)$ is complete for all nonanalytic algebra sequences $M$.
Every compact set which is a finite union of uniformly regular sets
in the sense of Dales and Davie is good.

{\bf Theorem 5:}
Let $X$ be a good compact plane set, let $\phi:X \rightarrow X$ and
let $\phi \in D^{\infty}(X)$.
Suppose that there is a point $b \in X$ such that $|\phi'(b)|>1$
and $\phi(b)$ has an external circular tangent.
Then there is a nonanalytic
sequence $M_n$ with $\phi \in D(X,M)$ such that $\phi$
does not induce an endomorphism of $D(X,M)$.

{\bf Proof:} Again we choose $M_n$ inductively, but we also choose
a sequence of complex numbers $c_n$ in the complement of $X$. 
For any $c$ off $X$ we define $F_c$ in $D^\infty(X)$ by 
$F_c(z)=1/(z-c)$.
We now begin with $M_0=M_1=1$ and any $c_0$, $c_1$ off $X$. 
For $n\geq 2$, and having chosen $M_0,\dots,M_{n-1}$
and $c_0, \cdots, c_{n-1},$
we choose $c_n$ and $M_n$ as follows. 
Choose $M_n$ such that (i) $M_n > \unorm {\phi^{(n)}} \infty /2^n$,
(ii) $M_n\geq (n!)^2$,
(iii) ${M_n\over{M_k M_{n-k}}} \geq {n \choose k}$ for $0 \leq k \leq n-1$,
and such that for all $k<n$, 
$$\unorm {F_{c_k}^{(n)}} \infty / 
M_n < 2^{-n} \unorm {F_{c_k}^{(k)}} \infty / M_k.$$

For the same $C_{n,\phi,m}$ as before,
we have, for all $F \in D^\infty(X)$ and all $z \in X$,
$$|(F \circ \phi)^{(n)}(z)| \geq |F^{(n)}(\phi(z))||\phi'(z)|^n
    -\sum_{m=0}^{n-1}{C_{n,\phi,m}\|F^{(m)}\|_\infty}.$$
We now consider the functions $F_c$ for suitable $c$. Since
$\phi(b)$ has an external circular tangent,
we can find $c$ in the complement of $X$ 
arbitrarily close to $\phi(b)$ and such that 
$|\phi(b)-c|=\min\{|z-c|: z \in X\}$. Choosing such $c=c_n$ 
close enough to $\phi(b)$ and considering the nature of the
derivatives of $F_c$, we can arrange that
$$|(F_c \circ \phi)^{(n)}(\phi(b))| \geq 
\frac{1}{2}\|F_c^{(n)}\|_\infty|\phi'(b)|^n$$
and also that $$\unorm {F_c^{(n)}} \infty / M_n \geq 
\sum_{m=0}^{n-1}{\unorm {F_c^{(m)}} \infty / M_m}.$$
The inductive choice may now proceed. 

     Having chosen the sequences 
$M_n$ and $c_n$, we see that $M_n$ is a nonanalytic  sequence, 
and that $\phi \in D(X,M)$. Also, for $n \geq 2$, 
$F_{c_n}$ is  in $D(X,M)$,
with $\unorm{F_{c_n}}{D(X,M)} \leq 3 \unorm {F_{c_n}^{(n)}} \infty / M_n $.
Further we have that $F_{c_n} \circ \phi \in D(X,M)$, so that
$$\unorm{F_{c_n} \circ \phi}{D(X,M)} \geq 
|(F_{c_n} \circ \phi)^{(n)}(\phi(b))| /M_n \geq 
\frac{1}{2 M_n}\|F_{c_n}^{(n)}\|_\infty|\phi'(b)|^n
\geq \frac{|\phi'(b)|^n}{6} \unorm{F_{c_n}}{D(X,M)}.$$
Since any endomorphism of $D(X,M)$ must be bounded, and $|\phi'(b)|>1$ 
it follows that $\phi$ does not induce an endomorphism of $D(X,M)$.


{\bf Part C}


    We now look at the special case where $X=\bar{\Delta}$.
We will use the property that elements in $D(\bar{\Delta},M)$
are analytic on $\Delta$ to obtain
nearly complete results for the case when the
self-maps 
are analytic on $\bar{\Delta}.$

     Suppose $\phi$ is analytic on $\Delta$ and
$\phi:\Delta \rightarrow \Delta$. It is well known that
if $\phi_n$ denotes the $n^{th}$
iterate of $\phi$, then unless $\phi$ is a rotation, 
there is a unique point $z_0 \in \bar{\Delta}$
such that $\phi_n (z) \rightarrow z_0$ for all $z \in \Delta.$
This point is known as the {\it Denjoy-Wolff point of $\phi.$} 
If $\phi$ is continuous at $z_0$, and certainly if $\phi$ is 
analytic on $\bar{\Delta}$, then the Denjoy-Wolff point, $z_0$, is
a fixed point of $\phi.$
If, in addition, $\phi'(z_0)$ exists, then $|\phi'(z_0)| \leq 1.$
Further, for all other fixed points $z'$ of $\phi$, $\phi'(z') > 1,$
whenever the derivative at $z'$ exists. Thus for nonanalytic
weights $M_n$, Theorem 4 shows that
if $\phi \in D(\bar{\Delta},M)$, $\phi$ is analytic on $\bar{\Delta}$,
 and $\phi:\bar{\Delta} \rightarrow \bar{\Delta}$ has two (or more) fixed points, then $\phi$ does not induce an
endomorphism of $D(\bar{\Delta},M).$ 
    
    We also recall that an inner function is an analytic function
$\phi$ on $\Delta$ such that $|\phi(z)| \leq 1$ and $|\phi(e^{i\theta})|=1$
for almost all $\theta.$

     The following classification of
analytic functions $\phi:\bar{\Delta} \rightarrow \bar{\Delta}$
in terms of fixed points was shown in \cite{Hp}.

{\bf Proposition:} Suppose $\phi:\bar{\Delta} \rightarrow \bar{\Delta}$
is analytic on $\bar{\Delta}.$ Then the following mutually
exclusive cases occur.
\begin{enumerate} 
   \item $\phi$ is inner.
    \item $\phi$ is not inner and for all integers $N$ there is
no fixed point of $\phi_N$ on the unit circle.
    \item There is a positive integer $N$ for which 
$S_N=\{\phi_N(w):|\phi_N(w)|=1\}$ is finite, nonempty and every
$z \in S_N$ is a fixed point of $\phi_N.$
 \begin{enumerate}
    \item $\phi$ has no fixed point in $\Delta$.
   \begin{enumerate}
    \item $\phi'(z') < 1$ for some $z' \in S_N$ and $\phi'_N(z) > 1$
for $z' \neq z \in S_N.$
    \item $\phi'(z')=1$ for some $z' \in S_N$ and $\phi'_N(z) > 1$
for $z' \neq z \in S_N.$
    \end{enumerate}
   \item $\phi_N$ (and hence $\phi$) has a fixed point in $\Delta$
and $\phi'_N(z) > 1$ for all $z \in S_N$.
\end{enumerate}
\end{enumerate}

{\bf Theorem 6:} Suppose $M_n$ is a nonanalytic weight sequence
and $\phi$ is an inner function in $D(\bar{\Delta},M).$ Then $\phi$
does not induce an endomorphism of $D(\bar{\Delta},M)$ unless
$\phi$ is a  constant or a rotation.

{\bf Proof:}

     Suppose that $\phi$ satisfies the hypothesis. Since $\phi$ is inner
and is continuous on $\bar{\Delta}$, it follows that $\phi$ must be
a finite Blaschke product. In particular, $\phi$ is
analytic on $\bar{\Delta}$.
We observe that if $\|\phi'\|_{\infty} > 1$, then there
exists $b \in \bar{\Delta} \setminus \Delta$ such that $|\phi(b)|=1$
and $|\phi'(b)| > 1.$ Theorem 4 shows that $\phi$ does not induce
an endomorphism of $D(\bar{\Delta},M)$ in this case. Thus we may assume
that
$\|\phi'\|_{\infty} \leq 1.$ Let $N(\phi)$ denote the number of
zeros of $\phi$ in $\Delta$. Then
\[\dis N(\phi)=\frac{1}{2\pi i}\int_{\bar{\Delta}\setminus \Delta}
\frac{\phi'(z)}{\phi(z)} dz \leq \frac{1}{2 \pi} \int_{\bar{\Delta}
\setminus \Delta} \|\phi'\|_{\infty} |dz| \leq 1\]
since $|\phi(z)|=1$ on the unit circle.
Therefore $N(\phi)=0$ or $1$. 
If $N(\phi)=0$, then $\phi$ is constant.
     Otherwise  $\phi$ is a  M\"{o}bius function,
so $\dis \phi(z)=e^{i\theta}\frac{z - \alpha}{1 - \bar{\alpha}z}$
for some real $\theta$ and $\alpha \in \Delta.$  Clearly,
\[\phi'(z)=e^{i\theta}\frac{1-|\alpha|^2}
{(1-\bar{\alpha}z)^2}.\]
     If $\alpha \neq 0,$ then $\dis
\phi'(\frac{\alpha}{|\alpha|})=\frac{1+|\alpha|}{1-|\alpha|} > 1.$
Theorem 4 then implies that the map
$\phi$ does not induce an endomorphism unless $\alpha =0$ in which
case $\phi$ has the form $\dis \phi(z)=e^{i\theta}z.$

   We remark that Theorem 6 also implies that the only automorphisms
of $D(\bar{\Delta},M)$ are those induced by rotations.

   Another consequence is to show that the completeness of the algebra
is needed  in Theorem 4. For, if we let  $D_r(\bar{\Delta},M)$ denote
the set of rational functions with poles off $\bar{\Delta}$, and $M_n$
a nonanalytic weight,  then the map $\dis T:Tf(z)=
f(\frac{2z-1}{z-2})$ induces an unbounded automorphism of the
incomplete normed algebra $D_r(\bar{\Delta},M).$ For, if $T$ were
bounded, then $T$ would extend to a bounded automorphism of $D(\bar{\Delta},M)$
induced by $\dis \phi(z) = \frac{2z-1}{z-2}$, which is impossible 
by Theorem 6.

{\bf Theorem 7:} Suppose $M_n$ is a nonanalytic weight,
$\phi$ is analytic on $\bar{\Delta}$,
$\phi:\bar{\Delta} \rightarrow \bar{\Delta}$, 
and $\phi$ is not an inner function.
Suppose that for all positive integers $N$ there is no fixed point
of $\phi_N$ on the unit circle. Then for some $N_1$, $\phi_{N}$ 
induces an endomorphism of $D(\bar{\Delta},M)$ for $N \geq N_1.$

{\bf Proof:}

     It follows from the hypothesis that $\phi$ has a (unique) fixed
point $z_0$ in $\Delta$. Since $|\phi'(z_0)| < 1$, a compactness argument
shows that $\phi_N(z) \in \Delta$ for all $z \in \bar{\Delta}$, $N$ large. 
Theorem 2 implies that
$\phi_N$ induces an endomorphism of $D(\bar{\Delta},M)$.

     However, we have the following example. Let 
$\dis \phi(z)=\frac{1-z^3}{2}.$ Here $\phi(-1)=1$ and
$\phi'(-1)=-\frac{3}{2}.$ Therefore $\phi$ does not induce an
endomorphism of $D(\bar{\Delta},M)$, while $\phi_2(z)=\phi(\phi(z))=
\frac{7-3z^3+3z^6-z^9}{16}$ does induce an endomorphism since
$\|\phi_2\|_{\infty} < 1.$

{\bf Theorem 8:} Suppose $\phi$ is analytic on $\bar{\Delta}$,
$\phi:\bar{\Delta} \rightarrow \bar{\Delta}$, 
$\phi \in  D(\bar{\Delta},M)$ and $\phi$ is not an inner function.
If $\phi$ has a fixed point in $\Delta$ and $\phi_N$ has a fixed
point on $\bar{\Delta} \setminus \Delta$ for some $N$, then
$\phi$ does not induce an endomorphism of $D(\bar{\Delta},M)$.

{\bf Proof:}

    Let $N$ be a positive integer and $z_1$ be a fixed point of
$\phi_N$ on $\bar{\Delta} \setminus \Delta.$ Then $\phi_N'(z_1) > 1$
and so $\phi_N$ does not induce an endomorphism of 
$D(\bar{\Delta},M).$ Clearly if $\phi_N$ does not induce an
endomorphism, then $\phi$ does not induce an endomorphism.

In the case when all of the fixed points of $\phi$ lie on
the boundary of $\bar{\Delta},$ we have the following.

{\bf Theorem 9:} Suppose $\phi$ is analytic on $\bar{\Delta}$,
$\phi \in D(\bar{\Delta},M)$, $\phi:\bar{\Delta}
\rightarrow \bar{\Delta}$ and $\phi$ has fixed points only on $\bar{\Delta}
\setminus \Delta.$ 

(i) If $\phi$ or $\phi_N$, for some $N$, has more than one fixed point 
on $\bar{\Delta} \setminus \Delta$,
then $\phi$ does not induce an endomorphism of $D(\bar{\Delta},M).$

(ii) If $\phi$ has exactly one fixed point $z_1$ on $\bar{\Delta} \setminus
\Delta$ and $\phi'(z_1) < 1$, then $\phi_N$ induces  an endomorphism of $D(\bar{\Delta},M)$ for all $N$ large enough.

{\bf Proof:}

(i) Let $z_0$ be the Denjoy-Wolff  point of $\phi_N$ in $\bar{\Delta} 
\setminus \Delta$. Then $\phi_N'(z_0) \leq 1.$ Then if $z_1$ is a
second fixed point of $\phi_N$, we have $|\phi_N(z_1)|=1$ and $|\phi_N'(z_1)| > 1.$
Theorem 4 then implies that $\phi_N$ and hence $\phi$
does not induce an endomorphism of
$D(\bar{\Delta},M).$

(ii) Say $\phi$ has exactly one fixed point $z_1$ on $\bar{\Delta} \setminus \Delta$ and $\phi'(z_1)  = r' < 1.$ 
Let $U$ be a neighborhood of $z_1$
for which $\dis |\phi'(z)| < r=\frac{1+r'}{2}$ and 
$\phi(z) \in U$ whenever $z \in U.$  
Let $A_n=\{z:\phi_n(z) \in U\}.$ 
By the Denjoy-Wolff Theorem, 
$\bigcup_{n=0}^{\infty} A_n=\bar{\Delta}$. 
Since the set $\{A_n\}$ is nested, the compactness of $\bar{\Delta}$ shows that
there exists $K_0$ so that $\phi_k(\bar{\Delta}) \subset U$ for
$k > K_0.$  Then for all $z \in \bar{\Delta}$,
and integers $N > K_0$,
\[\phi_N'(z)=\phi'(\phi_{N-1}(z)) \cdots \phi'(\phi_{K_0}(z))\phi'(\phi_{{K_0}-1}(z)) \cdots \phi'(\phi(z))\phi'(z)\]
so that \[\|\phi_N'\|_{\infty} < r^{N-K_0}(\|\phi'\|_{\infty})^{K_0}=
r^N(\frac{\|\phi'\|_{\infty}}{r})^{K_0}.\]
Then if $N > K_0$ is such that $\dis r^N(\frac{\|\phi'\|_{\infty}}{r})^{K_0} < 1$,
we have $\|\phi'_N\|_{\infty} < 1$ and so $\phi_N$ induces an
endomorphism by Theorem 1.

Remark: There is still one case which is not resolved, namely,
for $\phi \in D(\bar{\Delta},M)$, $\phi:\bar{\Delta} \rightarrow
\bar{\Delta}$ with $\phi(z_1)=z_1 \in \bar{\Delta} \setminus \Delta$, $\phi'(z_1)=1$, $
|\phi(z)|<1$
for $z \neq z_1$, but $\|\phi'\|_{\infty} > 1.$
An example of such $\phi$ is the following.
\[\phi(z)=\frac{1}{2}[z + \frac{(1+i)z-1}{z+(i-1)}].\]
   However, Theorem 3 shows that for sufficiently rapidly growing
$M_n$ this map $\phi$ induces an endomorphism of $D(\bar{\Delta},M).$

  We conclude with some open questions.

   1. What is the complete answer for the $\phi's$ which do not satisfy
the nice analytic properties that have been imposed?

   2. The automorphisms of $D([0,1],M)$ are induced by
$\phi(x)=x$ or $\phi(x)=1-x$ and the automorphisms of
$D(\bar{\Delta},M)$ are induced by rotations. Can the
automorphisms of other $D(X,M)$ be easily described?
We remark that there are perfect, compact sets $X$ such that the only
automorphism of $D(X,M)$ is the identity operator.

  3. What is the situation when $X$ is a disconnected set such as
the Cantor set?

  4. The spectra of composition operators on Banach spaces of analytic
functions on various domains have been studied in great detail.
Can those methods be used to determine the spectra of endomorphisms
of the algebras we have been considering?

\vspace{1.2in}
\vspace{.3in}

{\sf  School of Mathematical Sciences

 University of Nottingham 

 Nottingham NG7 2RD, England

 email: Joel.Feinstein@nottingham.ac.uk

and 

 Department of Mathematics

 University of Massachusetts at Boston 

 100 Morrissey Boulevard 

 Boston, MA 02125-3393

 email: hkamo@cs.umb.edu

\vspace{.5in}

This research was supported by EPSRC grant GR/M31132}

\end{document}